\documentclass[11pt]{article}

\usepackage[utf8]{inputenc}
\usepackage[T1]{fontenc}
\usepackage{amsmath,amssymb,amsthm,mathtools}
\usepackage{mathrsfs}
\usepackage[margin=1in]{geometry}
\usepackage[hidelinks,colorlinks=true,linkcolor=blue,citecolor=blue]{hyperref}
\usepackage{booktabs}
\usepackage{tikz}
\usepackage{xcolor}

\newtheorem{thm}{Theorem}[section]
\newtheorem{lem}[thm]{Lemma}

\newtheorem{corol}[thm]{Corollary}
\newtheorem{prop}[thm]{Proposition}
\theoremstyle{remark}

\theoremstyle{definition}

\newtheorem{problem}[thm]{Problem}
\newtheorem{defn}[thm]{Definition}
\newtheorem{example}[thm]{Example}

\newcommand{\sort}{\mathbf{sort}}

\newcommand{\keywords}[1]{%
  \par\smallskip
  \noindent\textbf{Keywords. }#1
}
\newcommand{\subjclass}[2][]{%
  \par\smallskip
  \noindent\textbf{2020 MSC. }#2
}

\title{\textbf{Iterating the Lehmer code on inversion sequences: Catalan fixed points and finite stabilization}}

\author{
Julian~Allagan$^{1}$,
Shanzhen~Gao$^{2}$,
Benjamin~Testart$^{3,*}$
}

\date{}

\begin{document}

\maketitle

\begin{center}
\small

$^{1}$Department of Mathematics,\\
University of Maryland Eastern Shore, Princess Anne, MD 21853, USA\\
\texttt{jallagan@umes.edu}

\vspace{4pt}

$^{2}$Department of Computer Information Systems,\\
Virginia State University, Petersburg, VA 23806, USA\\
\texttt{sgao@vsu.edu}

\vspace{4pt}

$^{3}$Universit\'e de Lorraine, CNRS, Inria, LORIA,\\
F-54000 Nancy, France\\
\texttt{benjamin.testart@loria.fr}

\vspace{6pt}

$^{*}$Corresponding author

\end{center}

\begin{abstract}
We study an operator $\Theta$ on finite integer sequences, where $\Theta(\sigma)_i$ counts the entries to the left of $\sigma_i$ that are strictly smaller than $\sigma_i$. This operator is a variant of the so-called Lehmer code. For every sequence $\sigma$, the image $\Theta(\sigma)$ is an inversion sequence, and the restriction of $\Theta$ to permutations of $[0,n-1]$ is a bijection onto inversion sequences of length $n$. We characterize the fixed points of $\Theta$ by avoidance of the pattern $101$ together with a saturation condition, prove that they are counted by the Catalan numbers, and give an explicit recursive bijection with Dyck paths. We also show that the sequences whose first $\Theta$-image is fixed are precisely those avoiding both $101$ and $201$. Finally, we prove finite stabilization for all inversion sequences, exhibit a family attaining the maximal stabilization time, and show that the second stabilization level is not closed under classical patterns.
\end{abstract}

\keywords{
inversion sequence,
Lehmer code,
permutation,
pattern avoidance,
Catalan numbers,
semi-Baxter numbers,
fixed point,
stabilization
}

\subjclass[2020]{
05A05,
05A15,
05A19
}

\section{Introduction}

Permutation patterns and their avoidance form a central part of enumerative combinatorics, with well-established connections to Catalan classes, sorting operators, lattice paths, logical descriptions of permutation classes, and related discrete structures~\cite{Albert_Bouvel_Feray_2020,Kitaev_2011,Simion_Schmidt_1985,Vatter_2015}. An inversion sequence of length $n$ is a sequence $\sigma=(\sigma_1,\ldots,\sigma_n)$ such that $0\leq \sigma_i<i$ for every $i\in[1,n]$. Inversion sequences are a well-known coding of permutations; a classical bijection from permutations onto reversed inversion sequences $(\sigma_n,\ldots,\sigma_1)$ is the Lehmer code \cite{Lehmer_1960}, which records the number of smaller entries on the right of each entry of a permutation.

This paper studies the increasing left-to-right counterpart of the Lehmer code\footnote{This is equivalent to reversing a permutation, applying the Lehmer code to it, then reversing the resulting sequence. Thus the fixed points and stabilization dynamics of our operator correspond to those of the Lehmer code under reversal. This convention ensures that images conform to the standard definition of an inversion sequence.}, and applies this coding to arbitrary integer sequences rather than restricting it to permutations. For a finite integer sequence ${\sigma=(\sigma_1,\ldots,\sigma_n)}$, define $\Theta(\sigma)_i=\#\{j<i:\sigma_j<\sigma_i\}$. Thus $\Theta(\sigma)_i$ counts earlier smaller entries\footnote{Note that in the literature, permutations are often turned into inversion sequences through another variant of the Lehmer code, counting earlier \emph{greater} entries; this operator has very different fixed points from the one we study (its fixed points are sequences of zeros).}, equivalently the occurrences of the classical permutation pattern $01$ ending at position $i$. For every finite integer sequence $\sigma$, the image $\Theta(\sigma)$ is an inversion sequence, since $0\leq\Theta(\sigma)_i<i$ for each $i$. When $\Theta$ is restricted to permutations of $[0,n-1]$, it gives a classical bijection onto inversion sequences of length $n$.

We analyze the dynamics of iterating $\Theta$ on permutations and inversion sequences. In this respect, $\Theta$ may be viewed as a stabilization operator analogous in spirit to classical sorting operators studied in permutation-pattern theory, such as the stack-sorting map \cite{Albert_Bouvel_Feray_2020,West_1993}, although the dynamics and fixed points considered here are fundamentally different. The first result characterizes the fixed points of $\Theta$ by two explicit conditions: avoidance of the pattern $101$ and a saturation condition on every value that appears. This characterization yields an enumeration of the fixed points of $\Theta$ according to their length by the Catalan numbers. We give an explicit recursive bijection between fixed points of length $n$ and Dyck paths of semilength $n$. We next prove that a sequence maps to a fixed point after one application of $\Theta$ precisely when it avoids both $101$ and $201$. Inversion sequences avoiding these two patterns are known to be enumerated by the semi-Baxter numbers~\cite{Bouvel_Guerrini_Rechnitzer_Rinaldi_2018,Martinez_Savage_2018}. 
Since permutations cannot contain the pattern $101$, the permutations whose first $\Theta$-image is fixed are the $201$-avoiding permutations, known to be counted by the Catalan numbers \cite{Knuth_1998, Simion_Schmidt_1985}. We prove finite stabilization under $\Theta$ for all inversion sequences with a sharp bound of $n-2$ iterations, show that the second stabilization level is not closed under classical patterns, and provide an upper bound for the number of permutations which stabilize within $k$ iterations.

\medskip

After releasing this work on arXiv, we were informed that repeated iterations of $\Theta$ were previously studied by \v Suni\'k~\cite{Sunik_2003}. Fortunately, most of our results are still original; there are only small overlaps with this previous work.
\begin{itemize}
    \item \v Suni\'k shows that all inversion sequences of length $n$ stabilize in $O(n^2)$ iterations of $\Theta$; we improve this to a tight upper bound of $n-2$ iterations in Propositions~\ref{prop:stabilization_criterion} and \ref{prop:sharp_stabilization}.
    \item \v Suni\'k characterizes the fixed points recursively through a tree-based construction (similar to a generating tree), then derives their Catalan enumeration from this construction. We provide a more direct characterization of fixed points in Theorem~\ref{thm:fixed_point_characterization}, and recover their Catalan enumeration independently, as well as providing a bijection with Dyck paths in Theorem~\ref{thm:fixed_dyck_bijection}.
\end{itemize}

\section{Preliminaries}

For integers $a\leq b$, write $[a,b]=\{i\in\mathbb Z:a\leq i\leq b\}$. For a finite sequence $\sigma=(\sigma_1,\ldots,\sigma_n)$, write $|\sigma|=n$. The empty sequence is allowed and has length $0$.

\begin{defn}
Let $S_n$ denote the set of permutations of $[0,n-1]$, written as sequences of length $n$. Let $I_n$ denote the set of inversion sequences of length $n$, that is, sequences $\sigma=(\sigma_1,\ldots,\sigma_n)$ satisfying $0\leq \sigma_i<i$ for every $i\in[1,n]$.
\end{defn}

\begin{defn}
Let $\sigma$ be a finite integer sequence. An occurrence of the pattern $101$ is a triple of indices $a<b<c$ such that $\sigma_a=\sigma_c>\sigma_b$. An occurrence of the pattern $201$ is a triple of indices $a<b<c$ such that $\sigma_a>\sigma_c>\sigma_b$. A sequence avoids a pattern if it has no occurrence of that pattern.
\end{defn}

\begin{defn}
For a finite integer sequence $\sigma$, define $\Theta(\sigma)$ by $\Theta(\sigma)_i=\#\{j<i:\sigma_j<\sigma_i\}$ for each $i\in[1,|\sigma|]$.
\end{defn}

\begin{example}
For $\sigma = (4,8,1,3,5,4)$, we have $\Theta(\sigma) = (0,1,0,1,3,2)$.
\end{example}

The following proposition describes well-known properties of $\Theta$; a proof is provided for the sake of completeness.
\begin{prop}\label{prop:theta_basic}
For every finite integer sequence $\sigma$, one has $\Theta(\sigma)\in I_{|\sigma|}$. Moreover, for every $n\geq0$, the restriction $\Theta:S_n\to I_n$ is a bijection.
\end{prop}

\begin{proof}
For each $i$, the value $\Theta(\sigma)_i$ counts a subset of $[1,i-1]$, so $0\leq \Theta(\sigma)_i<i$. Hence $\Theta(\sigma)\in I_{|\sigma|}$. It remains to prove the bijection on $S_n$. The case $n=0$ is immediate. Let $\tau\in I_n$. We reconstruct a unique permutation $\sigma\in S_n$ with $\Theta(\sigma)=\tau$ from right to left. Suppose $\sigma_{i+1},\ldots,\sigma_n$ have already been chosen, and let $A_i$ be the set of unused values. Since the values in $A_i$ are exactly the values that must occupy positions $1,\ldots,i$, the entry $\sigma_i$ must have precisely $\tau_i$ smaller values in $A_i$. Because $0\leq\tau_i<i=|A_i|$, there is a unique such choice, namely the $(\tau_i+1)$-st smallest element of $A_i$. This determines $\sigma_i$ uniquely at every step. The constructed sequence is a permutation of $[0,n-1]$, and by construction exactly $\tau_i$ entries among positions $1,\ldots,i-1$ are smaller than $\sigma_i$ for each $i$. Thus $\Theta(\sigma)=\tau$. Existence and uniqueness follow, so $\Theta:S_n\to I_n$ is bijective.
\end{proof}

\begin{defn}
For $n,k\geq0$, define:
\begin{enumerate}
    \item[(i)] $FP_n=\{\sigma\in I_n:\Theta(\sigma)=\sigma\}$, the set of
    fixed points of $\Theta$ of length $n$;
    \item[(ii)] $S_{n,k}=\{\sigma\in S_n:\Theta^k(\sigma)\in FP_n\}$, the set
    of permutations in $S_n$ whose $k$-th $\Theta$-image is fixed;
    \item[(iii)] $I_{n,k}=\{\sigma\in I_n:\Theta^k(\sigma)\in FP_n\}$, the set
    of inversion sequences in $I_n$ whose $k$-th $\Theta$-image is fixed.
\end{enumerate}
\end{defn}

Since every element of $FP_n$ is fixed by $\Theta$, we have $S_{n,k}\subseteq S_{n,k+1}$ and $I_{n,k}\subseteq I_{n,k+1}$ for all $n,k\geq0$. Proposition~\ref{prop:theta_basic} also gives $\Theta(S_{n,k})=I_{n,k-1}$ for every $k\geq1$: the inclusion $\Theta(S_{n,k})\subseteq I_{n,k-1}$ follows directly from the definitions, and the reverse inclusion follows by taking the unique preimage in $S_n$ under the bijection $\Theta:S_n\to I_n$.

\begin{prop}\label{prop:order_preserving}
Let $\sigma$ be an integer sequence of length $n$, and let $\tau = \Theta(\sigma)$. For all $p<q \in [1,n]$, the following hold: if $\sigma_p\leq\sigma_q$, then $\tau_p\leq\tau_q$; if $\sigma_p<\sigma_q$, then $\tau_p<\tau_q$.
\end{prop}

\begin{proof}
If $\sigma_p\leq\sigma_q$, then every index $i<p$ with $\sigma_i<\sigma_p$ also satisfies $i<q$ and $\sigma_i<\sigma_q$. Hence the set counted by $\tau_p$ is contained in the set counted by $\tau_q$, so $\tau_p\leq\tau_q$. If $\sigma_p<\sigma_q$, the same containment holds, and the index $p$ gives an additional contribution to $\tau_q$ that is not counted in $\tau_p$. Therefore $\tau_q\geq\tau_p+1$, so $\tau_p<\tau_q$.
\end{proof}

Proposition~\ref{prop:order_preserving} above shows that $\Theta$ preserves the relations $\leq$ and $<$ between two entries read from left to right. Note however this does not hold for the relations $=$, $\geq$, or $>$. For instance, if $\sigma = (2,0,1,2)$, then $\tau = \Theta(\sigma) = (0,0,1,2)$; we have $\sigma_1 > \sigma_3$ and $\sigma_1 = \sigma_4$, yet $\tau_1 < \tau_3$ and $\tau_1 < \tau_4$.

\section{Basic properties of \texorpdfstring{$\Theta$}{Θ} on inversion sequences}

Throughout this section, let $n\geq0$, let $\sigma\in I_n$, and $\tau=\Theta(\sigma)$. All statements involving a position $p$ are understood with $p\in[1,n]$.

\begin{prop}\label{prop:zeros_preserved}
If $\sigma_p=0$, then $\tau_p=0$.
\end{prop}

\begin{proof}
By definition, $\tau_p=\#\{i<p:\sigma_i<\sigma_p\}$. Since $\sigma_p=0$ and every entry of an inversion sequence is nonnegative, the set being counted is empty.
\end{proof}

\begin{prop}\label{prop:value_increase}
For every $p\in[1,n]$, one has $\sigma_p\leq\tau_p$.
\end{prop}

\begin{proof}
Since $\sigma\in I_n$, we have $\sigma_p<p$. For each $i\in[1,\sigma_p]$, the inversion-sequence condition gives $\sigma_i<i\leq\sigma_p$, and hence $\sigma_i<\sigma_p$. Thus the $\sigma_p$ indices $1,\ldots,\sigma_p$ all contribute to $\tau_p$, so $\tau_p\geq\sigma_p$.
\end{proof}

\begin{prop}\label{prop:local_stability}
If $\sigma_p=\tau_p$, then $\Theta(\tau)_p=\tau_p$.
\end{prop}

\begin{proof}
Let $i<p$. If $\sigma_i<\sigma_p$, then Proposition~\ref{prop:order_preserving} gives $\tau_i<\tau_p$. If $\sigma_i\geq\sigma_p$, then Proposition~\ref{prop:value_increase} gives $\tau_i\geq\sigma_i\geq\sigma_p=\tau_p$. Hence $\{i<p:\sigma_i<\sigma_p\}=\{i<p:\tau_i<\tau_p\}$. The two sets have the same cardinality, so $\Theta(\tau)_p=\tau_p$.
\end{proof}

\begin{prop}\label{prop:stabilization_criterion}
Let $p\in[1,n]$ and $k\geq1$. If $\Theta^k(\sigma)_p\leq k$, then $\Theta^{k-1}(\sigma)_p=\Theta^k(\sigma)_p$.
\end{prop}

\begin{proof}
Set $\sigma^{(r)}=\Theta^r(\sigma)$. We prove the claim by induction on $k$. For $k=1$, Proposition~\ref{prop:value_increase} gives $\sigma_p\leq\sigma^{(1)}_p\leq1$. If $\sigma_p=0$, then Proposition~\ref{prop:zeros_preserved} gives $\sigma^{(1)}_p=0$. If $\sigma_p=1$, then the preceding inequality forces $\sigma^{(1)}_p=1$. Thus $\sigma_p=\sigma^{(1)}_p$.

Assume now that $k\geq2$ and that the result has been proved for $k-1$. Proposition~\ref{prop:value_increase}, applied to $\sigma^{(k-1)}$, gives $\sigma^{(k-1)}_p\leq\sigma^{(k)}_p$. If the inequality were strict, then $\sigma^{(k-1)}_p<\sigma^{(k)}_p\leq k$, so $\sigma^{(k-1)}_p\leq k-1$. By the induction hypothesis, $\sigma^{(k-2)}_p=\sigma^{(k-1)}_p$. Proposition~\ref{prop:local_stability}, applied to $\sigma^{(k-2)}$ in place of $\sigma$, then gives $\sigma^{(k)}_p=\sigma^{(k-1)}_p$, a contradiction. Hence $\sigma^{(k-1)}_p=\sigma^{(k)}_p$.
\end{proof}

\begin{prop}\label{prop:finite_stabilization}
For every $n\geq2$ and every $\sigma\in I_n$, the sequence $\Theta^{n-2}(\sigma)$ is fixed by $\Theta$.
\end{prop}

\begin{proof}
For each $p\in[1,n]$, the entry $\Theta^{n-1}(\sigma)_p$ belongs to $[0,p-1]$, and hence is at most $n-1$. Applying Proposition~\ref{prop:stabilization_criterion} with $k=n-1$ gives $\Theta^{n-2}(\sigma)_p=\Theta^{n-1}(\sigma)_p$. Since this holds for every position $p$, we have $\Theta^{n-2}(\sigma)=\Theta^{n-1}(\sigma)$, so $\Theta^{n-2}(\sigma)$ is fixed by $\Theta$.
\end{proof}

Two immediate consequences follow: $I_{n,k}=I_n$ for every $k\geq n-2$ and $n\geq2$, and, by the bijectivity of $\Theta:S_n\to I_n$, one also has $S_{n,k}=S_n$ for every $k\geq n-1$ and $n\geq2$; moreover, Proposition~\ref{prop:sharp_stabilization} shows this is sharp, since $S_{n,n-2}\neq S_n$ for every $n\geq3$.

Example~\ref{ex:extremal_orbit} works out a representative member of an extremal
family; Proposition~\ref{prop:sharp_stabilization} shows that
$\sigma=(0,1,2,\ldots,n-3,0,1)$ attains the bound of
Proposition~\ref{prop:finite_stabilization} for every $n\ge3$.

\begin{example}\label{ex:extremal_orbit}
Take $n=5$ and $\sigma=(0,1,2,0,1)$. The first four entries form the fixed
point $(0,1,2,0)$, so under iteration only the final coordinate moves:
\[
(0,1,2,0,1)\ \xrightarrow{\ \Theta\ }\ (0,1,2,0,2)
\ \xrightarrow{\ \Theta\ }\ (0,1,2,0,3)
\ \xrightarrow{\ \Theta\ }\ (0,1,2,0,4),
\]
and $(0,1,2,0,4)$ is fixed. At each step the final entry $v$ is replaced by the
number of earlier entries smaller than $v$; because the value $0$ occurs twice
among the first four entries, this count equals $v+1$ rather than $v$, so the
final coordinate strictly increases until it reaches $n-1$. Stabilization therefore
occurs after $n-2=3$ iterations, the maximum permitted by
Proposition~\ref{prop:finite_stabilization}.
\end{example}

We now show that the stabilization bound from Proposition~\ref{prop:finite_stabilization} is attained.
\begin{prop}\label{prop:sharp_stabilization}
For every $n\geq3$, the inversion sequence $\sigma=(0,1,2,\ldots,n-3,0,1)$ satisfies $\Theta^{\,n-3}(\sigma)\neq\Theta^{\,n-2}(\sigma)$. Consequently, the bound of Proposition~\ref{prop:finite_stabilization} is sharp.
\end{prop}

\begin{proof}
We prove by induction on $k$ that $\Theta^k(\sigma)=(0,1,2,\ldots,n-3,0,k+1)$ for every $0\leq k\leq n-2$.

The assertion is immediate for $k=0$. Assume it holds for some $k<n-2$. The first $n-1$ entries of $\Theta^k(\sigma)$ form the fixed point $(0,1,2,\ldots,n-3,0)$, so only the final coordinate can change under another application of $\Theta$. The last entry of $\Theta^k(\sigma)$ is $k+1$. Among the preceding entries, exactly $k+2$ positions have values smaller than $k+1$: the positions $1,\ldots,k+1$, whose values are $0,1,\ldots,k$, together with the additional occurrence of $0$ at position $n-1$. Hence the final coordinate of $\Theta^{k+1}(\sigma)$ is $k+2$. Therefore $\Theta^{k+1}(\sigma)=(0,1,2,\ldots,n-3,0,k+2)$, completing the induction.

Taking $k=n-3$ gives $\Theta^{\,n-3}(\sigma)=(0,1,2,\ldots,n-3,0,n-2)$, whereas $\Theta^{\,n-2}(\sigma)=(0,1,2,\ldots,n-3,0,n-1)$. These sequences are distinct. Since Proposition~\ref{prop:finite_stabilization} gives stabilization by step $n-2$, this example attains the maximal possible stabilization time.
\end{proof}

\section{Fixed points}\label{sec:fixedpoints}

\begin{defn}
Let $\sigma\in I_n$. An entry $\sigma_i$ is called \emph{saturated} if $\sigma_i=i-1$.
\end{defn}

\begin{prop}\label{prop:image_leftmax_saturated}
Let $\sigma$ be a nonempty integer sequence, let $p$ be the position of the leftmost maximum of $\sigma$, and let $\tau=\Theta(\sigma)$. Then $p$ is the position of the rightmost saturated entry of $\tau$.
\end{prop}

\begin{proof}
Since $p$ is the position of the leftmost maximum of $\sigma$, every entry to the left of $p$ is strictly smaller than $\sigma_p$. Hence every index $1,\ldots,p-1$ contributes to $\tau_p$, and therefore $\tau_p=p-1$. Thus the entry at position $p$ is saturated.

Now let $j>p$. Since $\sigma_p\geq\sigma_j$, the index $p$ does not contribute to $\tau_j$. Hence at most the $j-2$ indices in $[1,j-1]\setminus\{p\}$ can contribute to $\tau_j$, so $\tau_j\leq j-2<j-1$. Therefore no saturated entry of $\tau$ occurs to the right of $p$, proving that $p$ is the rightmost saturated position of $\tau$.
\end{proof}

\begin{prop}\label{prop:preimage_characterization}
Let $\sigma$ and $\sigma'$ be two integer sequences of length $n \geq 0$. Then $\Theta(\sigma)=\Theta(\sigma')$ if and only if,
for every pair $p<q$, one has $\sigma_p<\sigma_q$ if and only if
$\sigma'_p<\sigma'_q$.
\end{prop}

\begin{proof}
Suppose first that $\sigma$ and $\sigma'$ have the same increasing pairs. Then, for each position $q$, the subsets $\{p<q:\sigma_p<\sigma_q\}$ and $\{p<q:\sigma'_p<\sigma'_q\}$ of $[1,q-1]$ are equal. Hence they have the same cardinality, and therefore $\Theta(\sigma)=\Theta(\sigma')$.

Conversely, assume $\Theta(\sigma)=\Theta(\sigma')=\tau$. We prove by induction on $n$ that $\sigma$ and $\sigma'$ have the same increasing pairs. The assertion is immediate for $n\leq1$, so let $n\geq2$.

Let $p$ be the position of the leftmost maximum of $\sigma$. By
Proposition~\ref{prop:image_leftmax_saturated}, $p$ is the rightmost saturated
position of $\tau$. Applying the same proposition to $\sigma'$ shows that $p$ is
also the position of the leftmost maximum of $\sigma'$. We first compare all pairs involving position $p$. If $i<p$, then $\sigma_i<\sigma_p$ and $\sigma'_i<\sigma'_p$, since $p$ is the leftmost maximum in both sequences. If $j>p$, then $\sigma_j\leq\sigma_p$ and $\sigma'_j\leq\sigma'_p$, since $p$ is a maximum in both sequences. Thus, for every $r\neq p$, the comparison between positions $p$ and $r$ has the same truth value in $\sigma$ as in $\sigma'$. Delete position $p$ from $\sigma,\sigma'$, and $\tau$, and denote the resulting sequences by $\widehat{\sigma},\widehat{\sigma}'$, and $\widehat{\tau}$, respectively. We claim that $\Theta(\widehat{\sigma})=\widehat{\tau}=\Theta(\widehat{\sigma}')$.

We verify the $\Theta$-relations after deletion. For a position originally
to the left of $p$, neither the entry nor any earlier entry is affected by the
deletion, so the corresponding $\Theta$-coordinate is unchanged. For a position
$q>p$, the deleted entry $\sigma_p$ satisfies $\sigma_p\geq\sigma_q$ and hence
does not contribute to $\Theta(\sigma)_q$; removing it therefore deletes no index
from the set counted by $\Theta(\sigma)_q$. After reindexing, the corresponding
coordinate of $\Theta(\widehat{\sigma})$ equals $\tau_q$. Hence
$\Theta(\widehat{\sigma})=\widehat{\tau}$, and in particular
$\widehat{\tau}\in I_{n-1}$ by Proposition~\ref{prop:theta_basic}. The identical
argument applied to $\sigma'$ gives $\Theta(\widehat{\sigma}')=\widehat{\tau}$.

By the induction hypothesis, $\widehat{\sigma}$ and $\widehat{\sigma}'$ have the
same increasing pairs. Together with the comparisons involving position $p$
verified above, this shows that $\sigma_a<\sigma_b$ if and only if
$\sigma'_a<\sigma'_b$ for every pair $a<b$ of original positions. Hence $\sigma$
and $\sigma'$ have the same increasing pairs, completing the proof.
\end{proof}

\begin{corol}\label{corol:fixed_preimage}
Let $\sigma$ be an integer sequence, and let $\tau=\Theta(\sigma)$. Then $\tau$ is fixed by $\Theta$ if and only if, for every $p<q$, one has $\sigma_p<\sigma_q$ if and only if $\tau_p<\tau_q$.
\end{corol}

\begin{proof}
Since $\Theta(\sigma)=\tau$, Proposition~\ref{prop:preimage_characterization} gives $\Theta(\tau)=\tau$ if and only if $\sigma$ and $\tau$ have exactly the same increasing pairs.
\end{proof}

\begin{thm}\label{thm:fixed_point_characterization}
Let $\sigma\in I_n$. Then $\sigma$ is fixed by $\Theta$ if and only if the following two conditions hold: \textup{(i)} $\sigma$ avoids $101$; \textup{(ii)} whenever a value $v$ appears in $\sigma$, one has $\sigma_{v+1}=v$.
\end{thm}

\begin{proof}
Let $\tau=\Theta(\sigma)$.

Assume first that $\sigma$ is fixed. If $\sigma$ contained an occurrence of $101$ at positions $a<b<c$, then $\sigma_a=\sigma_c>\sigma_b$. Every index contributing to $\tau_a$ would also contribute to $\tau_c$, and the index $b$ would contribute to $\tau_c$ but not to $\tau_a$. Hence $\tau_c>\tau_a$, contradicting $\tau=\sigma$ and $\sigma_a=\sigma_c$. Thus $\sigma$ avoids $101$.

Now let $v$ be a value appearing in $\sigma$, and choose $p$ with $\sigma_p=v$. By Proposition~\ref{prop:value_increase}, $\tau_p\geq v$. Since $\tau_p=\sigma_p=v$, equality holds. Note that the indices $1,\ldots,v$ already contribute $v$ units to $\tau_p$, hence no index $r$ with $v<r<p$ can satisfy $\sigma_r<v$. If $p=v+1$, then $\sigma_{v+1}=v$. Otherwise, $v+1<p$, so taking $r=v+1$ in the preceding sentence gives $\sigma_{v+1}\geq v$. Since $\sigma\in I_n$, one also has $\sigma_{v+1}<v+1$, and therefore $\sigma_{v+1}=v$.

Conversely, assume that $\sigma$ satisfies \textup{(i)} and \textup{(ii)}. Fix $p$ and write $v=\sigma_p$. By Proposition~\ref{prop:value_increase}, $\tau_p\geq v$. Condition \textup{(ii)} gives $\sigma_{v+1}=v$. If $r$ satisfies $v+1<r<p$ and $\sigma_r<v$, then the indices $v+1<r<p$ form an occurrence of $101$, since $\sigma_{v+1}=\sigma_p=v>\sigma_r$, contradicting \textup{(i)}. Hence every index between $v+1$ and $p$ has value at least $v$. However, the indices $1,\ldots,v$ all have values smaller than $v$, by the inversion-sequence condition. We conclude that exactly $v$ entries of $\sigma$ contribute to $\tau_p$, therefore $\tau_p=v=\sigma_p$. Since $p$ was arbitrary, $\Theta(\sigma)=\sigma$.
\end{proof}

\begin{lem}\label{lem:fixed_decomposition}
For $n\geq1$, every $\sigma\in FP_n$ has a unique decomposition determined by $v=\sigma_n$ of the form $\sigma=(\alpha,\beta_1+v,\ldots,\beta_{n-v}+v)$, where $\alpha\in FP_v$, $\beta\in FP_{n-v}$, and $\beta_{n-v}=0$. Conversely, every such pair $(\alpha,\beta)$ produces an element of $FP_n$.
\end{lem}

\begin{proof}
Let $\sigma\in FP_n$ and set $v=\sigma_n$. By Theorem~\ref{thm:fixed_point_characterization}, $\sigma_{v+1}=v$. Since $\sigma$ avoids $101$, no index $r$ with $v+1<r<n$ can satisfy $\sigma_r<v$. Thus every entry in positions $v+1,\ldots,n$ is at least $v$. Define $\alpha=(\sigma_1,\ldots,\sigma_v)$ and $\beta=(\sigma_{v+1}-v,\ldots,\sigma_n-v)$. Then $\beta_{n-v}=0$.

The sequence $\alpha$ belongs to $I_v$, and it satisfies the two fixed-point conditions of Theorem~\ref{thm:fixed_point_characterization}; hence $\alpha\in FP_v$. Likewise, $\beta\in I_{n-v}$. Any occurrence of $101$ in $\beta$ would give an occurrence of $101$ in $\sigma$, and if a value $w$ appears in $\beta$, then the value $v+w$ appears in $\sigma$. The fixed-point condition for $\sigma$ gives $\sigma_{v+w+1}=v+w$, which is exactly $\beta_{w+1}=w$. Hence $\beta\in FP_{n-v}$.

Conversely, let $\alpha\in FP_v$ and let $\beta\in FP_{n-v}$ with $\beta_{n-v}=0$. Define $\sigma=(\alpha,\beta_1+v,\ldots,\beta_{n-v}+v)$. This is an inversion sequence because $\alpha\in I_v$ and, for each $j\in[1,n-v]$, one has $\beta_j+v\leq j-1+v$. The first block has all entries smaller than $v$, while the second block has all entries at least $v$, so no occurrence of $101$ can use equal entries from both blocks. Since both blocks separately avoid $101$, the sequence $\sigma$ avoids $101$.

It remains to check the saturation condition. If a value $u<v$ appears in $\sigma$, it appears in $\alpha$, and $\alpha_{u+1}=u$ gives $\sigma_{u+1}=u$. If a value $u=v+w$ appears in the second block, then $w$ appears in $\beta$, and $\beta_{w+1}=w$ gives $\sigma_{u+1}=\sigma_{v+w+1}=v+w=u$. Hence $\sigma$ satisfies both conditions of Theorem~\ref{thm:fixed_point_characterization}, so $\sigma\in FP_n$. The value $v=\sigma_n$ uniquely determines the split, proving uniqueness.
\end{proof}

\begin{lem}\label{lem:fixed_ending_zero}
For $n\geq1$, the map $\rho\mapsto(\rho,0)$ is a bijection from $FP_{n-1}$ onto the set of elements of $FP_n$ whose last entry is $0$.
\end{lem}

\begin{proof}
If $\rho\in FP_{n-1}$, then $(\rho,0)$ is an inversion sequence. It avoids $101$, because the final entry is $0$ and cannot be the final term of a $101$ occurrence. The saturation condition is unchanged for all values appearing in $\rho$, and the value $0$ is already saturated at position $1$. Hence $(\rho,0)\in FP_n$.

Conversely, if $\sigma\in FP_n$ and $\sigma_n=0$, then deleting the last entry preserves avoidance of $101$ and preserves the saturation condition for every value that remains. Therefore $(\sigma_1,\ldots,\sigma_{n-1})\in FP_{n-1}$. The two constructions are inverse to each other.
\end{proof}

\begin{thm}\label{thm:catalan_enumeration}
For every $n\geq0$, one has $|FP_n|=C_n=\frac{1}{n+1}\binom{2n}{n}$.
\end{thm}

\begin{proof}
Let $A_n=|FP_n|$. Clearly $A_0=1$. For $n\geq1$, decompose each $\sigma\in FP_n$ according to Lemma~\ref{lem:fixed_decomposition}, with $v=\sigma_n$. For a fixed $v\in[0,n-1]$, the first block may be chosen in $A_v$ ways. The second block is an element of $FP_{n-v}$ ending in $0$, and Lemma~\ref{lem:fixed_ending_zero} shows that these are counted by $A_{n-v-1}$. Hence $A_n=\sum_{v=0}^{n-1}A_vA_{n-1-v}$. This is the Catalan recurrence with initial value $A_0=1$. Therefore $A_n=C_n=\frac{1}{n+1}\binom{2n}{n}$ for every $n\geq0$.
\end{proof}

\begin{figure}[ht]
\centering
\begin{tikzpicture}
    \begin{scope}[scale = 0.9]
    \draw[color = gray] \foreach \x in {0,...,3} {(0,\x) -- (18,\x)} \foreach \x in {0,...,18} {(\x,0) -- (\x,3)};
    \draw[line width = 2pt] (0,0) -- (1,1);
    \draw[line width = 2pt, color = red] (1,1) -- (2,2) -- (3,3) -- (4,2) -- (5,1) -- (6,2) -- (7,1);
    \draw[line width = 2pt] (7,1) -- (8,0);
    \draw[line width = 2pt, color = blue] (8,0) -- (9,1) -- (10,0) -- (11,1) -- (12,2) -- (13,3) -- (14,2) -- (15,3) -- (16,2) -- (17,1) -- (18,0);

    \node at (4, 3.5) {\huge \textcolor{red}{$A$}};
    \node at (13, 3.5) {\huge \textcolor{blue}{$B$}};
    \end{scope}
\end{tikzpicture}
\caption{The Dyck path obtained by applying the function $\Phi$ of Theorem~\ref{thm:fixed_dyck_bijection} to the sequence $(0,1,1,3,3,5,6,3,3) \in FP_9$, and its decomposition into $A = \Phi(0,1,1)$ and $B = \Phi(0,0,2,3,0)$.}
\label{fig:Dyck_path}
\end{figure}
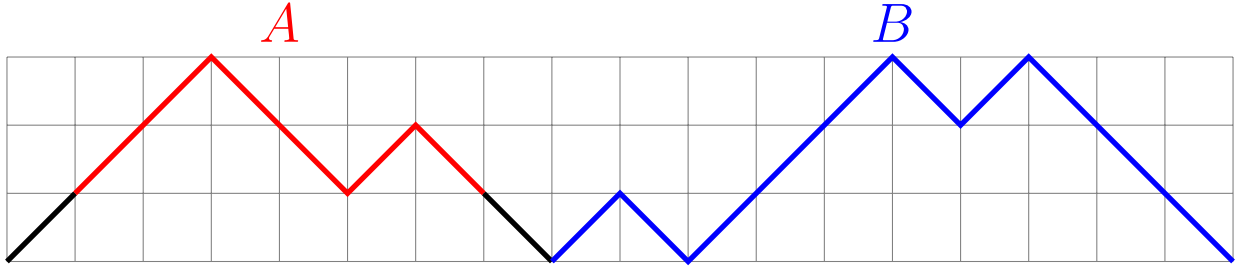

The Catalan enumeration of $FP_n$ can be strengthened to an explicit recursive bijection with Dyck paths. Figure~\ref{fig:Dyck_path} pictures this bijection on an example.

\begin{thm}[A recursive Catalan bijection]\label{thm:fixed_dyck_bijection}
For every $n\geq0$, there is a bijection $\Phi:FP_n\to\mathcal D_n$, where $\mathcal D_n$ denotes the set of Dyck paths of semilength $n$.
\end{thm}

\begin{proof}
We proceed by induction on $n$. For $n=0$, we define $\Phi(\emptyset)=\emptyset$.
Now let $n\geq 1$, and assume that, for every $m<n$, the map $\Phi:FP_m \to \mathcal D_m$ is a bijection. We define $\Phi$ on $FP_n$ and show that it is a bijection.

Let $\sigma\in FP_n$, and write $v=\sigma_n$.
By Lemma~\ref{lem:fixed_decomposition}, $\sigma$ admits a unique decomposition into $\sigma = (\alpha,\beta_1+v, \dots,\beta_{n-v}+v)$, where $\alpha\in FP_v$, $\beta\in FP_{n-v}$, and $\beta_{n-v}=0$.
Since $v<n$, the induction hypothesis applies to $\alpha$. Moreover, by
Lemma~\ref{lem:fixed_ending_zero}, deleting the final zero from $\beta$ gives a unique element $\gamma\in FP_{n-v-1}$.
Since $n-v-1<n$, the induction hypothesis also applies to $\gamma$.
We therefore define $\Phi(\sigma) = U\Phi(\alpha)D\Phi(\gamma)$.
The path $\Phi(\alpha)$ has semilength $v$, while $\Phi(\gamma)$ has semilength $n-v-1$. Hence $\Phi(\sigma)$ has semilength $n$, and thus $\Phi(\sigma)\in\mathcal D_n$.

It remains to prove that $\Phi$ is a bijection. Let $P\in\mathcal D_n$ be
nonempty. By the unique first-return decomposition of a Dyck path, there are
unique $v\in\{0,\ldots,n-1\}$, $A\in\mathcal D_v$, and $B\in\mathcal D_{n-v-1}$ such that $P=UADB$.
Since $v<n$ and $n-v-1<n$, the induction hypothesis implies that there are
unique $\alpha\in FP_v$ and $\gamma\in FP_{n-v-1}$ such that $\Phi(\alpha)=A$ and $\Phi(\gamma)=B$.
By Lemma~\ref{lem:fixed_ending_zero}, appending a final zero to $\gamma$ gives an element $\beta=(\gamma,0)\in FP_{n-v}$.
Let $\sigma = (\alpha,\beta_1+v, \dots,\beta_{n-v}+v)$.
By Lemma~\ref{lem:fixed_decomposition}, we have $\sigma\in FP_n$. By the
definition of $\Phi$, we have $\Phi(\sigma) = U\Phi(\alpha)D\Phi(\gamma) = UADB = P$.
Thus $\Phi$ is surjective. Since both $FP_n$ and $\mathcal D_n$ are counted by $C_n$, we conclude that $\Phi$ is a bijection.
\end{proof}

Theorem~\ref{thm:fixed_dyck_bijection} places the fixed points of $\Theta$ within the classical Catalan family of Dyck paths; further Catalan interpretations may be found in Stanley~\cite{Stanley_1999}.

\section{One-step stabilization}\label{sec:stabilization}

The characterization of fixed points obtained in Section~\ref{sec:fixedpoints} yields a complete description of the first stabilization level through pattern avoidance.

\begin{thm}\label{thm:avoid_101_201}
Let $\sigma$ be an integer sequence of length $n$. Then $\Theta(\sigma) \in FP_n$ if and only if $\sigma$ avoids both patterns $101$ and $201$.
\end{thm}

\begin{proof}
Let $\tau=\Theta(\sigma)$. Assume first that $\sigma$ contains the pattern $101$ or $201$ (or both). Let $a < b < c$ be a triple of indices such that $\sigma_a \geq \sigma_c > \sigma_b$, i.e. $(\sigma_a, \sigma_b, \sigma_c)$ is an occurrence of 101 or 201, and such that $a$ is minimal. By the minimality of $a$, we have $\sigma_i < \sigma_c$ for every $i < a$. Therefore every index $i<a$ contributes to $\tau_c$, and the index $b$ contributes as well because $\sigma_b<\sigma_c$. Consequently $\tau_c\geq a$. On the other hand, $\tau_a\leq a-1$ by Proposition~\ref{prop:theta_basic}. It follows that $\tau_a<\tau_c$, while $\sigma_a\geq\sigma_c$. Hence the equivalence in Corollary~\ref{corol:fixed_preimage} fails for the pair $(a,c)$, and therefore $\tau$ is not fixed. We have proved by contrapositive that if $\tau$ is fixed, then $\sigma$ avoids both $101$ and $201$.

Conversely, assume that $\sigma$ avoids both patterns. Let $p<q$.
\newline
\textit{Case 1.} Suppose $\sigma_p<\sigma_q$. By Proposition~\ref{prop:order_preserving}, one has $\tau_p<\tau_q$.
\newline
\textit{Case 2.} Suppose $\sigma_p\geq\sigma_q$. We claim that no index $i$ with $p<i<q$ satisfies $\sigma_i<\sigma_q$. Indeed, if such an index existed, then $(p,i,q)$ would form an occurrence of $101$ when $\sigma_p=\sigma_q$, and an occurrence of $201$ when $\sigma_p>\sigma_q$. Both possibilities contradict the hypothesis. Consequently, every index contributing to $\tau_q$ also contributes to $\tau_p$, and therefore $\tau_p\geq\tau_q$. Combining the two cases yields $\sigma_p<\sigma_q$ if and only if $\tau_p<\tau_q$. Since this holds for every pair $p<q$, Corollary~\ref{corol:fixed_preimage} implies that $\tau$ is fixed.
\end{proof}

\begin{corol}
Inversion sequences of length $n$ whose $\Theta$-image is fixed under $\Theta$ are counted by $|I_{n,1}|$, with $|I_{0,1}|=|I_{1,1}|=1$ and, for $n\geq2$,
\[
|I_{n,1}|=\frac{11n^2+11n-6}{(n+4)(n+3)}|I_{n-1,1}|+\frac{(n-3)(n-2)}{(n+4)(n+3)}|I_{n-2,1}|.
\]
\end{corol}
\begin{proof}
    By Theorem~\ref{thm:avoid_101_201}, the set $I_{n,1}$ consists precisely of the inversion sequences avoiding both $101$ and $201$. Martinez and Savage~\cite{Martinez_Savage_2018} showed that this class of inversion sequences is enumerated by the semi-Baxter numbers (entry A117106 of \cite{OEIS}). These numbers were studied in depth by Bouvel, Guerrini, Rechnitzer, and Rinaldi~\cite{Bouvel_Guerrini_Rechnitzer_Rinaldi_2018}, who found the D-finite recurrence formula above.
\end{proof}

Corollary~\ref{corol:Catalan_permutations} below can be seen as an immediate consequence of Theorem~\ref{thm:catalan_enumeration}, since $\Theta$ is a bijection from $S_{n,1}$ onto $FP_n$. This result can also be obtained through Theorem~\ref{thm:avoid_101_201}, which yields additional insight on the structure of permutations in $S_{n,1}$.

\begin{corol}\label{corol:Catalan_permutations}
Permutations whose $\Theta$-image is fixed under $\Theta$ are counted by the Catalan numbers: for every $n\geq0$, $|S_{n,1}|=C_n=\frac{1}{n+1}\binom{2n}{n}$.
\end{corol}
\begin{proof}
Since all permutations avoid the pattern 101, Theorem~\ref{thm:avoid_101_201} implies that the set $S_{n,1}$ consists precisely of the $201$-avoiding permutations in $S_n$. Permutations in $S_n$ avoiding a 3-letter permutation pattern are well-known to be counted by $C_n$, see for instance \cite{Simion_Schmidt_1985}. \qedhere
\end{proof}

\section{Stabilization dynamics}

Theorem~\ref{thm:avoid_101_201} determines the first stabilization level. For $k \geq 2$, the classes $S_{n,k}$ are less transparent. Although membership in $S_{n,k}$ is defined by iterating $\Theta$, the resulting conditions do not reduce naturally to avoidance of a finite set of classical permutation patterns. The case $k=2$ already reflects this obstruction.

\begin{prop}\label{prop:S2_not_classical}
The union $\bigcup_{n\geq0}S_{n,2}$ is not closed under taking classical permutation patterns. Consequently, it is not a classical permutation class.
\end{prop}

\begin{proof}
Let $\rho=(2,3,0,1)\in S_4$. Then $\Theta(\rho)=(0,1,0,1)$, $\Theta^2(\rho)=(0,1,0,2)$, and $\Theta^3(\rho)=(0,1,0,3)$. Since $\Theta^2(\rho)$ is not fixed, $\rho\notin S_{4,2}$.

Now let $\pi=(3,4,1,0,2)\in S_5$. One has $\Theta(\pi)=(0,1,0,0,2)$ and $\Theta^2(\pi)=(0,1,0,0,4)$. The latter sequence is fixed by $\Theta$, so $\pi\in S_{5,2}$. However, the subsequence of $\pi$ in positions $1,2,4,5$ is $(3,4,0,2)$, whose reduction is $(2,3,0,1)=\rho$. Thus an element of $S_{5,2}$ contains a classical pattern that does not belong to $S_{4,2}$. Therefore $\bigcup_{n\geq0}S_{n,2}$ is not closed under classical patterns.
\end{proof}

By Theorem~\ref{thm:avoid_101_201}, a permutation $\pi$ belongs to $S_{n,2}$ if and only if $\Theta(\pi)$ avoids both $101$ and $201$. Translating this condition directly into the entries of $\pi$ yields the following explicit, though nonlocal, characterization of $S_{n,2}$.

\begin{prop}\label{prop:Sn2}
A permutation $\pi\in S_n$ does not belong to $S_{n,2}$ if and only if
there exist indices $\ell<r$ such that 
$\#\{\,i<\ell:\pi_r<\pi_i<\pi_\ell\,\}
\ge
\#\{\,m:\ell<m<r,\ \pi_m<\pi_r\,\}
\ge 1$.
\end{prop}

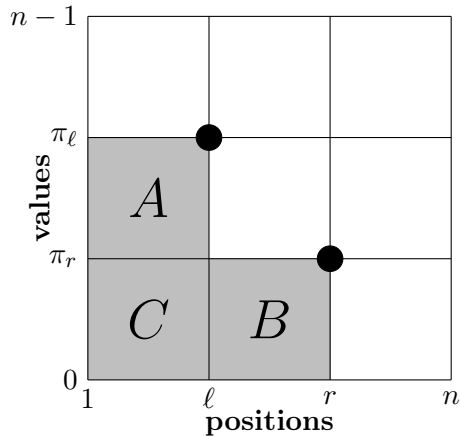
\begin{figure}[ht]
\centering
\begin{tikzpicture}
\begin{scope}[scale=0.8]
    \fill[color=lightgray] (0,2) -- (0,4) -- (2,4) -- (2,2) -- (0,2);
    \node at (1, 3) {\huge $A$};
    \fill[color=lightgray](2,0) -- (2,2) -- (4,2) -- (4,0) -- (2,0);
    \node at (3, 1) {\huge $B$};
    \fill[color=lightgray] (0,0) -- (0,2) -- (2,2) -- (2,0) -- (0,0);
    \node at (1, 1) {\huge $C$};
    \filldraw (2, 4) circle(6pt);
    \filldraw (4, 2) circle(6pt);

    \node[left] at (0,0) {$0$};
    \node[left] at (0,2) {$\pi_r$};
    \node[left] at (0,4) {$\pi_\ell$};
    \node[left] at (0,6) {$n-1$};
    
    \draw (0,0) -- (6,0);
    \draw (0,2) -- (6,2);
    \draw (0,4) -- (6,4);
    \draw (0,6) -- (6,6);
    \draw (0,0) -- (0,6);
    \draw (2,0) -- (2,6);
    \draw (4,0) -- (4,6);
    \draw (6,0) -- (6,6);
    \node[rotate = 90] at (-0.8,3) {\textbf{values}};
    
    \node at (0,-0.3) {1};
    \node at (2,-0.3) {$\ell$};
    \node at (4,-0.3) {$r$};
    \node at (6,-0.3) {$n$};
    \node at (3,-0.7) {\textbf{positions}};
\end{scope}
\end{tikzpicture}
\caption{Illustration of Proposition~\ref{prop:Sn2} showing the regions counted by $A$, $B$, and $C$ in the permutation $\pi$.}
\label{fig:Sn2}
\end{figure}

\begin{proof}
Let $\tau=\Theta(\pi)$. Assume first that $\pi\notin S_{n,2}$. By
Theorem~\ref{thm:avoid_101_201}, the inversion sequence $\tau$
contains an occurrence of either $101$ or $201$. Hence there exist indices $\ell<q<r$ such that $\tau_\ell\ge\tau_r>\tau_q$.
Since Proposition~\ref{prop:order_preserving} implies that
$\pi_\ell<\pi_r$ would force $\tau_\ell<\tau_r$, one must have
$\pi_\ell>\pi_r$. As pictured on Figure~\ref{fig:Sn2}, define
\begin{align*}
A&=\#\{\,i<\ell:\pi_r<\pi_i<\pi_\ell\,\},\\
B&=\#\{\,m:\ell<m<r,\ \pi_m<\pi_r\,\},\\
C&=\#\{\,i<\ell:\pi_i<\pi_r\,\}.
\end{align*}

Because $\pi_\ell>\pi_r$, every index counted by $\tau_\ell$
either contributes to $C$ or contributes to $A$. Hence
$\tau_\ell=A+C$. Likewise, every index counted by $\tau_r$
either contributes to $C$ or contributes to $B$. Hence
$\tau_r=B+C$. Since $\tau_\ell\ge\tau_r$, it follows that
$A\ge B$. It remains to show that $B\ge1$.

If $\pi_q<\pi_r$, then $q$ itself contributes to $B$ since
$\ell<q<r$. Suppose instead that $\pi_q>\pi_r$. Every index $i<q$ satisfying $\pi_i<\pi_r$ also satisfies $\pi_i<\pi_q$ and
therefore contributes to $\tau_q$. Since $\tau_r>\tau_q$,
not all indices contributing to $\tau_r$ can lie before $q$.
Hence there exists an index $m$ with $q<m<r$ and
$\pi_m<\pi_r$. Such an index contributes to $B$. Thus $B\ge1$, and therefore $A\ge B\ge1$.

Conversely, suppose there exist indices $\ell<r$ such that
\[
A:=\#\{\,i<\ell:\pi_r<\pi_i<\pi_\ell\,\} \ge B:=\#\{\,m:\ell<m<r,\ \pi_m<\pi_r\,\}\ge1.
\]
Since $A\ge1$, there exists an index $i<\ell$ satisfying
$\pi_r<\pi_i<\pi_\ell$. In particular, $\pi_\ell>\pi_r$. Let
$ C=\#\{\,i<\ell:\pi_i<\pi_r\,\}.$
As before,
$\tau_\ell=A+C$
and
$\tau_r=B+C$.
Hence $\tau_\ell\ge\tau_r$. Since $B\ge1$, choose an index $m$ such that
$\ell<m<r$ and $\pi_m<\pi_r$.
By Proposition~\ref{prop:order_preserving},
$\tau_m<\tau_r$.
Consequently, $\tau_\ell\ge\tau_r>\tau_m$. Therefore $(\ell,m,r)$ forms an occurrence of either $101$ or $201$
in $\tau$. By Theorem~\ref{thm:avoid_101_201},
$\tau\notin I_{n,1}$.
Since $\tau=\Theta(\pi)$, it follows that
$\pi\notin S_{n,2}$.
\end{proof}

\subsection{Exact-round data}

For $\pi\in S_n$, define its stabilization time by $T_n(\pi)=\min\{k\geq0:\Theta^k(\pi)\in FP_n\}$. Thus $\pi\in S_{n,k}$ if and only if $T_n(\pi)\leq k$. Further, for $n\geq1$ and $k\geq0$, let $E_{n,k}=|S_{n,k}\setminus S_{n,k-1}|$, where $S_{n,-1}=\emptyset$. Equivalently, $E_{n,k}=|\{\pi\in S_n:T_n(\pi)=k\}|$, so $E_{n,k}$ counts the permutations that stabilize for the first time after exactly $k$ applications of $\Theta$. The numbers $E_{n,k}$ are now entry A399547 of \cite{OEIS}.

Table~\ref{tabl:stabilization_time} records the values of $E_{n,k}$ for
small $n$. These values were obtained by exhaustive enumeration of $S_n$
and direct iteration of $\Theta$.

\begin{table}[ht]
\begin{center}
\small
\begin{tabular}{c|rrrrrrrrrr}
$n\backslash k$ & 0 & 1 & 2 & 3 & 4 & 5 & 6 & 7 & 8 & 9 \\
\hline
1 & 1 \\
2 & 1 & 1 \\
3 & 1 & 4 & 1 \\
4 & 1 & 13 & 9 & 1 \\
5 & 1 & 41 & 62 & 15 & 1 \\
6 & 1 & 131 & 398 & 167 & 22 & 1 \\
7 & 1 & 428 & 2529 & 1706 & 343 & 32 & 1 \\
8 & 1 & 1429 & 16304 & 17079 & 4781 & 679 & 46 & 1 \\
9 & 1 & 4861 & 107795 & 172118 & 64385 & 12340 & 1314 & 65 & 1 \\
10 & 1 & 16795 & 733930 & 1768042 & 866136 & 210604 & 30720 & 2481 & 90 & 1
\end{tabular}
\end{center}
\caption{Values of $E_{n,k}$ for $0\leq k<n\leq10$.}
\label{tabl:stabilization_time}
\end{table}

The first column reflects that the identity permutation is the unique
fixed point in $S_n$, while the second column equals $C_n-1$, in
agreement with Corollary~\ref{corol:Catalan_permutations}. The final
nonzero entry in each row occurs at $k=n-1$. Indeed,
Proposition~\ref{prop:finite_stabilization} gives stabilization of every
inversion sequence by step $n-2$, so by the bijection
$\Theta:S_n\to I_n$ of Proposition~\ref{prop:theta_basic} every
permutation stabilizes by step $n-1$; the unique preimage of the extremal
inversion sequence in Proposition~\ref{prop:sharp_stabilization} attains
this bound, stabilizing for the first time at step $n-1$. Since
$E_{n,k}=|\{\pi\in S_n:T_n(\pi)=k\}|$, $|S_n|=n!$, and $T_n\leq n-1$, for a
uniformly random permutation in $S_n$ we have
$\mathbb{P}(T_n=k)=E_{n,k}/n!$. Hence
$\mathbb{E}[T_n]=\frac{1}{n!}\sum_{k=0}^{n-1}kE_{n,k}$. The values in
Table~\ref{tabl:stabilization_time} already indicate that typical
stabilization occurs well before the worst-case bound $n-1$; its
asymptotic behavior is considered in Problem~\ref{prob:average_growth}.

\subsection{Upper bound for \texorpdfstring{$|S_{n,k}|$}{|Snk|}}
While it is generally difficult to compute the exact value of $|S_{n,k}|$ for $k \geq 2$, the present section provides an upper bound which holds for any $k$.

For any integer sequence $\sigma$, let $\sort(\sigma)$ denote the sequence obtained by sorting the entries of $\sigma$ in increasing order.

\begin{lem} \label{lem_injectivity}
    Let $\sigma$ and $\sigma'$ be two integer sequences of length $n \geq 0$. If $\sort(\sigma) = \sort(\sigma')$ and $\Theta(\sigma) = \Theta(\sigma')$, then $\sigma = \sigma'$.
\end{lem}

\begin{proof}
We proceed by induction on $n$. For $n=0$, there is only the empty sequence, so the result is immediate.

Now assume that $n>0$, that the result holds for sequences of length $n-1$, and that $\sigma$ and $\sigma'$ satisfy $\sort(\sigma) = \sort(\sigma')$ and $\Theta(\sigma) = \Theta(\sigma')$. In particular, $\sigma$ and $\sigma'$ have the same multiset of values, by the hypothesis $\sort(\sigma) = \sort(\sigma')$. Since $\Theta(\sigma)_n$ counts the entries of $\sigma$ with value less than $\sigma_n$, the equality $\Theta(\sigma)_n = \Theta(\sigma')_n$ now implies that $\sigma_n = \sigma'_n$.

Let $\bar{\sigma}=(\sigma_1,\ldots,\sigma_{n-1})$ and $\bar{\sigma}'=(\sigma'_1,\ldots,\sigma'_{n-1})$ be the sequences obtained by deleting the final term from $\sigma$ and $\sigma'$. Since $\sigma_n=\sigma'_n$, we have $\sort(\bar{\sigma})=\sort(\bar{\sigma}')$.
Furthermore, for every $i<n$, the value $\Theta(\sigma)_i$ depends only on the terms $\sigma_j$ with $j<i$, hence $\Theta(\bar{\sigma})=\Theta(\bar{\sigma}')$.
By the induction hypothesis, it follows that $\bar{\sigma}=\bar{\sigma}'$. Together with $\sigma_n=\sigma'_n$, we conclude that $\sigma=\sigma'$.
\end{proof}

\begin{thm}\label{thm:Catalan_bound}
For all $n,k\geq0$, we have $|S_{n,k}|\leq (C_n)^k$, where $C_n=\frac{1}{n+1}\binom{2n}{n}$.
\end{thm}

\begin{proof}
    Fix $n$ and proceed by induction on $k$. For $k \leq 1$, we have $|S_{n,0}| = 1$ since only the permutation $(0,1, \dots, n-1)$ is fixed under $\Theta$, and $|S_{n,1}| = C_n$ by Corollary~\ref{corol:Catalan_permutations}.
    
    Assume now that $k \geq 2$, and that $|S_{n,k-1}| \leq (C_n)^{k-1}$. By Proposition~\ref{prop:theta_basic}, the identity ${|S_{n,k}| = |I_{n,k-1}|}$ holds. By Lemma~\ref{lem_injectivity}, we have $|I_{n,k-1}| = \#\{(\sort(\sigma), \Theta(\sigma))  :  \sigma \in I_{n,k-1}\}$. It follows that
    \begin{equation} \label{eq:sort_Theta}
        |S_{n,k}| \leq |\sort(I_{n,k-1})| \cdot |\Theta(I_{n,k-1})|.
    \end{equation}
    Note that any $\sigma \in I_n$ satisfies $0 \leq \sigma_j < j \leq i$ for all $j \leq i$, so at least $i$ entries of $\sigma$ are strictly smaller than $i$; hence $\sort(\sigma)_i < i$ for all $i \in [1,n]$, and $\sort(\sigma) \in I_n$. Consequently, the elements of $\sort(I_{n,k-1})$ are nondecreasing inversion sequences. The set of all nondecreasing inversion sequences of length $n$ is known to be counted by $C_n$, see for instance \cite[Section 2.14.1]{Martinez_Savage_2018} or \cite[Section 2.1]{Testart_2025}. Therefore $|\sort(I_{n,k-1})| \leq C_n$.
    Furthermore, ${\Theta(I_{n,k-1}) \subseteq I_{n,k-2} = \Theta(S_{n,k-1})}$, hence $|\Theta(I_{n,k-1})| \leq |S_{n,k-1}|$ by Proposition~\ref{prop:theta_basic}. Replacing both terms on the right-hand side of \eqref{eq:sort_Theta}, we obtain ${|S_{n,k}| \leq C_n \cdot |S_{n,k-1}|}$. By the induction hypothesis, we conclude that $|S_{n,k}| \leq (C_n)^k$.
\end{proof}

An interesting consequence of Theorem~\ref{thm:Catalan_bound} is that, for any fixed $k$, the numbers $|S_{n,k}|$ grow at most exponentially in $n$. Indeed, since $C_n \leq 4^n$ for all $n \geq 0$, Theorem~\ref{thm:Catalan_bound} implies that $|S_{n,k}| \leq (4^k)^n$.

\section{Conclusion and open problems}

We studied the action of the operator $\Theta$ on inversion sequences and permutations. We obtained a complete characterization of the fixed points of $\Theta$ in terms of avoidance of $101$ together with a saturation condition, established their Catalan enumeration, and constructed an explicit bijection with Dyck paths. We further characterized the inversion sequences that stabilize within one application of $\Theta$ as those avoiding both $101$ and $201$, recovering a semi-Baxter class, and showed that the corresponding permutations are precisely the $201$-avoiding permutations. We proved finite stabilization for all inversion sequences, established a sharp bound of $n-2$ iterations, and exhibited extremal sequences attaining this bound. Finally, we showed an upper bound of $\left (\frac{1}{n+1}\binom{2n}{n} \right )^k$ for the number of permutations of length $n$ which stabilize within $k$ iterations.

The results reveal a nontrivial interaction between inversion-sequence patterns, Catalan structures, and iterative dynamics. At the same time, they leave several natural questions open. The exact behavior of higher stabilization levels remains largely unexplored, both from enumerative and structural perspectives.

\begin{problem}\label{prob:average_growth}
Determine the asymptotic growth of $\mathbb{E}[T_n]$ for a uniformly
random permutation in $S_n$.
\end{problem}

Monte Carlo experiments for $n$ up to approximately $2000$ suggest a
power-law exponent close to $3/4$. The experiments also indicate
comparable growth of the mean and median of $T_n$. Determining the
limiting behavior of the distribution of $T_n$, under an appropriate
normalization, remains open.

\begin{problem}\label{prob:higher_generating}
For fixed $k\geq2$, find a way to efficiently compute the numbers $|I_{n,k}|$, or determine their asymptotic behavior.
\end{problem}

While the cases $k=0$ and $k=1$ admit Catalan and semi-Baxter descriptions, respectively, no comparable enumeration is presently known for higher stabilization levels.

\begin{problem}\label{prob:higher_intrinsic}
For fixed $k\geq2$, find structural properties characterizing the elements of $S_{n,k}$ or $I_{n,k}$.
\end{problem}

Proposition~\ref{prop:S2_not_classical} shows that the class $\bigcup_{n\geq0}S_{n,2}$ is not closed under classical pattern containment. Consequently, any satisfactory description of higher stabilization levels will likely require tools beyond the traditional theory of permutation classes.

More broadly, it would be interesting to determine whether the stabilization process admits a deeper probabilistic, geometric, or algebraic interpretation analogous to those arising from classical sorting operators and pattern-avoidance dynamics.

\section*{Acknowledgments}
The idea of repeatedly applying the Lehmer code was proposed during the pre-conference workshop of Permutation Patterns 2025. Participants of the workshop found the links to Catalan numbers and 201-avoiding permutations (cf Corollary~\ref{corol:Catalan_permutations}). We gratefully acknowledge their initial contributions, which motivated our further work. We also thank Krishna Menon and Emil Verkama for bringing the reference \cite{Sunik_2003} to our attention.

\end{document}